\documentclass[11pt]{article}

\usepackage{latexsym}
\usepackage{amssymb}

\newcommand{\R}{\mathbb R}
\newcommand{\h}{\mathcal H}
\newcommand{\C}{\mathbb C}
\newcommand{\la}{\lambda}

\newcommand{\n}{\mathfrak n}
\newcommand{\T}{\mathfrak T}

\renewcommand{\ll}{\lesssim}

\title{The Hilbert transform, rearrangements, \\
and logarithmic determinants}
\author{Vladimir Matsaev\thanks
{Supported by the Israel Science Foundation of the Israel Academy
of Sciences and Humanities under Grant No. 37/00-1.}, \ Iossif Ostrovskii,
\ and Mikhail Sodin$^{\ast}$}

\begin{document}

\date{}

\maketitle

This is an extended version of notes 
prepared for the talk at the conference 
``Rajchman-Zygmund-Marcinkiewicz 2000''. 
They are based on recent papers \cite{MOS} and \cite{MS2}
(see also \cite{MS1} and \cite{MS3}). 
The authors thank Professor Zelazko for the invitation 
to participate in the proceedings of this conference.

\bigskip\par\noindent{\bf \S 1}

\medskip\par\noindent
Let $g$ be a bounded measurable real-valued function on $\R$
with a compact support.

We shall use the following notations:

\smallskip\par\noindent The Hilbert transform of $g$:
$$
(\h g) (\xi) = \frac{1}{\pi} \int'_\R \frac{g(t)}{t-\xi}\, dt\,,
$$
the prime means that the integral is understood in the principal value 
sense at the point $t=\xi$.
\smallskip\par\noindent The (signed) distribution function of $g$:
$$
N_g(s)=
\left\{ \begin{array}{ll}
\ \  {\rm meas}\,\{x:\, g(x)> s\}\,, &  \mbox{if \ $s>0$}; \\ \\  
-{\rm meas}\,\{x:\, g(x)< s\}\,, &  \mbox{if \ $s<0$}.
\end{array}
\right.
$$ 

\smallskip\par\noindent The (signed) decreasing rearrangement of $g$:
$g_d$ is defined as the distribution function of $N_g$:
$g_d = N_{N_g}$. 

Less formally, the functions $N_g$ and $g_d$ can be also defined by the 
following properties: they are non-negative and non-increasing for $s>0$, 
non-positive and non-increasing for $s<0$, and 
$$
\int_\R \Phi (g(t))\, dt 
=\int_\R \Phi (s)\, dN_g(s)
=\int_\R \Phi(g_d(t))\, dt\,,
$$
for any function $\Phi$ such that at least one of the three integrals is 
absolutely convergent.

We shall use notation $A \lesssim B$, when $A\le C\cdot B$ for a positive 
numerical constant $C$. We shall write $A \lesssim_\lambda B$, if $C$ in 
the previous inequality depends on the parameter $\lambda $ only. 

\medskip\par\noindent{\bf Theorem~1.1 }{\em Let 
$g$ be a bounded measurable real-valued function with a compact
support. Then}
$$
||\h g_d||_{L^1} \le 4 ||\h g||_{L^1}\,.
\leqno (1.2)
$$

\smallskip\par\noindent
Hereafter, $L^1$ always means $L^1(\R)$.

\medskip\par\noindent{\bf Remarks: } 

\par\noindent{\bf 1.3 } Estimate (1.2) can be extended to a wider 
class of functions after an additional regularization of the Hilbert 
transform $\h g_d$ (see $\S3$ below). 

\par\noindent{\bf 1.4 } Probably, the constant $4$ on the RHS is not 
sharp. However, Davis' discussion in \cite{Davis} suggests that (1.2) 
ceases to hold without this factor on the RHS.

\par\noindent{\bf 1.5 } Theorem~1.1 yields a result of Tsereteli 
\cite{Tse} and Davis \cite{Davis}: {\em if $g\in \mbox{Re} H^1$, then 
$g_d$ is also in} $\mbox{Re}H^1$, {\em and} $||\h g_d||_{L^1} \lesssim 
||g||_{\mbox{Re}H^1}$, where $\mbox{Re} H^1$ is the real Hardy space on 
$\R$.  

\par\noindent{\bf 1.6 } Theorem~1.1 can be extended to functions defined 
on the unit circle $\mathbb T$. Let $g(t)$ be a bounded function on 
$\mathbb T$, $g_d$ be its signed decreasing rearrangement, and $\tilde g$ 
be a function conjugate to $g$:
$$
\tilde g(t) = \frac{1}{2\pi} \int'_{\mathbb T} 
g(\xi) \cot \frac{t-\xi}{2}\, d\xi\,.
$$
Then
$$
||\widetilde{g_d}||_{L^1(\mathbb T)} 
\le 4 ||\tilde{g}||_{L^1(\mathbb T)}\,. 
\leqno (1.7)
$$
Juxtapose this estimate with Baernstein's inequality \cite{Ba}:
$$
||\tilde{g}||_{L^1(\mathbb T)} 
\le ||\widetilde{g_s}||_{L^1(\mathbb T)}\,, 
\leqno (1.8)
$$
where $g_s$ is the symmetric decreasing rearrangement of $g$. 
In particular, if $g_s$ has a conjugate in $L^1$, then 
{\em any} rearrangement of $g$ has a conjugate in $L^1$, and if 
{\em some} rearrangement of $g$ has a conjugate in $L^1$, then 
the conjugate of $g_d$ is in $L^1$. 
We are not aware of a counterpart of Baernstein's inequality for the 
Hilbert transform and the $L^1(\R)$-norm.

\bigskip\par\noindent{\bf \S 2}

\medskip\par\noindent
Here, we shall prove Theorem~1.1. WLOG, we assume that
$$
\int_{\R} g(t)\, dt = 0\,,
\leqno (2.1)
$$
otherwise 
$$
(\h g)(\xi) = - \frac{1}{\pi \xi} \int_\R g(t)\,dt + O(1/\xi^2)\,, 
\qquad \xi\to\infty\,,
$$
and the $L^1$-norm on the RHS of (1.2) is infinite.

\medskip\par\noindent{\bf The first reduction:} instead of 
(1.2), we shall prove inequality
$$
||\h N_g||_{L^1} \le 2 ||\h g||_{L^1}\,,
\leqno (2.2)
$$
then its iteration gives (1.2). 

We introduce a (regularized) logarithmic determinant of $g$:
$$
u_g(z) \stackrel{def}= \int_\R K(z g(t))\, dt\,,
\qquad K(z)=\log|1-z|+\mbox{Re}(z)\,.
$$
This function is subharmonic in $\C$ and harmonic outside 
of $\R$. 

\medskip\par\noindent{\bf List of properties of $u_g$:}

\smallskip\par\noindent Since $g$ is a bounded function with a compact 
support, 
$$
u_g(z) = O(|z|^2)\,, \qquad z\to 0\,,
\leqno (2.3a)
$$
and by (2.1)
$$
u_g(z) = \int_\R \log|1-zg(t)|\, dt = O(\log |z|)\,, \qquad z\to\infty\,.
\leqno (2.3b)
$$
In particular,
$$
\int_\R \frac{|u_g(x)|}{x^2} < \infty\,.
\leqno (2.3c)
$$

\smallskip\par\noindent Next,
$$
\int_\R \frac{u_g(x)}{x^2}\,dx = 0\,.
\leqno (2.4)
$$
This follows from the Poisson representation:
$$
u_g(iy) = \frac{y}{\pi} \int_\R \frac{u_g(x)}{x^2+y^2}\, dy\,,
\qquad y>0\,.
$$
Dividing by $y$, letting $y\to 0$, and using (2.3a), we get (2.4).

\smallskip\par\noindent Further,
$$
u_g(1/t)=-\pi(\h N_g)(t)\,.
\leqno (2.5)
$$
Indeed, integrating by parts and changing variables, we obtain
for real $x$'s:
\begin{eqnarray*}
u_g(x) 
&=& \int_\R \log|1-xs|\, dN_g(s) \\
&=& x \int'_\R \frac{N_g(s)}{1-xs}\, ds \\ \\
&=& -\pi (\h N_g)(1/x)\,.
\end{eqnarray*}

\medskip
We have done the {\bf second reduction: }
Instead of (2.2), we shall prove inequality
$$
\int_\R \frac{u_g^-(x)}{x^2}\, dx \le 
\pi ||\h g||_{L^1}\,.
\leqno (2.6)
$$
Then combining (2.4) and (2.6), we get (2.2).

Now, we set
$$
f(t) = g(t) + i(\h g)(t)\,.
$$
This function has an analytic continuation into the upper 
half-plane:
$$
f(z) = \frac{1}{\pi i} \int_{\R} \frac{g(t)}{t-z}\, dt\,.
$$
We define the regularized logarithmic determinant of $f$ by the equation
$$
u_f(z) = \int_\R K(zf(t))\, dt\,.
\leqno (2.7)
$$
The {\em positivity} of this subharmonic function is central in our 
argument:

\par\noindent{\bf Lemma~2.8 (cf. \cite{Essen}) }
$$
u_f(z) \ge 0\,,
\qquad z\in\C\,.
$$

\smallskip\par\noindent{\em Proof of Lemma~2.8: } It suffices to consider 
such $z$'s that all solutions of the equation $zf(w)=1$ are simple and 
not real. Then 
\begin{eqnarray*}
u_f(z) &=& \mbox{Re}
\left\{
\int_\R \big[ \log(1-zf(t)) +zf(t) \big] dt 
\right\} \\
&=& \mbox{Re}
\left\{ z^2 
\int_\R \frac{tf(t)f'(t)}{1-zf(t)}\, dt
\right\} \\ 
&=& \mbox{Re} 
\left\{
2\pi i z^2 \sum_{\{w:\, zf(w)=1\}} \mbox{Res}_w
\left( \frac{\zeta f(\zeta) f'(\zeta)}{1-zf(\zeta)} \right)
\right\} \\
&=& 2\pi \sum_{\{w:\, zf(w)=1\}} \mbox{Im}(w) \ge 0\,.
\end{eqnarray*}
The application of the Cauchy theorem is justified since 
$f(\zeta)=O(1/\zeta^2)$ when $\zeta\to\infty$, $\mbox{Im}(\zeta)\ge 0$. 
Done.

\smallskip To complete the proof of the theorem, we shall use an argument 
borrowed from the perturbation theory of compact operators \cite{GK}. 
We use auxiliary functions $f_1=g+i|\h g|$ and
$$
u_1(z) = \int_\R \log\left|
\frac{1-zg(t)}{1-zf_1(t)}
\right|\, dt\,.
$$
Then on the real axis
$$
u_g(x) = u_1(x) + u_f(x)\,, \qquad x\in\R\,,
$$
so that $u_g(x)\ge u_1(x)$, or $u_g^-(x)\le u_1^-(x) = -u_1(x)$, since 
$u_1(x)\le 0$, $x\in \R$.

Next, we need an elementary inequality: if $w_1$, $w_2$ 
are complex numbers such that $\mbox{Re}(w_1) = \mbox{Re}(w_2)$ and
$|\mbox{Im}(w_1)|\le \mbox{Im}(w_2)$, then for all $z$ in the upper 
half-plane,
$$
\left| \frac{1-zw_1}{1-zw_2} \right| <1\,.
$$
Due to this inequality the function $u_1$ is non-positive in the upper 
half-plane. Since this function is harmonic in the upper 
half-plane, we obtain
\begin{eqnarray*}
\int_\R \frac{u_g^-(x)}{x^2}\, dx &\le& 
- \int_\R \frac{u_1(x)}{x^2}\, dx \\
&=& - \lim_{y\to 0}\int_\R \frac{u_1(x)}{x^2+y^2}\, dx \\ 
&\le& - \pi \lim_{y\to 0} \frac{u_1(iy)}{y} \\
&=& -\pi \lim_{y\to 0} \frac{1}{y} \int_\R \log
\left| \frac{1-iyg(t)}{1-iyg(t)+y|(\h g)(t)|}
\right|\, dt \\ 
&=&\pi \int_\R |(\h g)(t)|\, dt\,.
\end{eqnarray*}
This proves (2.6) and therefore the theorem. $\Box$

\bigskip\par\noindent{\bf \S 3}

\medskip\par\noindent  
Here, we will formulate a fairly complete version of estimate (2.2). The 
proof given in \cite{MS2} follows similar lines as above, however is 
essentially more technical. 

Now, we start with a real-valued measure $d\eta$ of finite variation on 
$\R$, and denote by $g=\h \eta$ its Hilbert transform. 
By $||\eta||$ we denote the 
total variation of the measure $d\eta$ on $\R$.   
Let 
$R_g=\h^{-1} N_g$ be a regularized inverse Hilbert transform of $N_g$:
$$
R_g(t) \stackrel{def}=
\lim_{\epsilon\to 0}\frac{1}{\pi} \int'_{|s|>\epsilon}
\frac{N_g(s)}{t-s}\, ds\,.
$$ 
The integral converges at infinity due to the Kolmogorov weak
$L^1$-type estimate 
$$
N_g(s) \lesssim ||\eta||/s\,, 
\qquad 0<s<\infty\,.
$$ 
Existence of the limit when $\epsilon\to 0$ (and $t\ne 0$) follows from
the Titchmarsh formula \cite{T} (cf. \cite{MS2}):
$$
\lim_{s\to 0} sN_g(s) = \frac{\eta (\R)}{\pi}\,.
$$

\medskip\par\noindent{\bf Theorem 3.1 } {\em Let $d\eta$ be a real measure 
supported by $\R$. Then
$$
\int_\R R_g^+(t)dt \le ||\eta_{\rm{a.c.}}||\,,
\leqno (3.2)
$$
$$
\int_\R R_g^-(t)dt \le ||\eta|| - |\eta(\R)|\,,
\leqno (3.3)
$$
and} 
$$
\int_\R R_g(t)dt =  |\eta (\R)| - ||\eta_{\rm{sing}}||\,.
\leqno (3.4)
$$

\medskip\par\noindent
{\bf Corollary 3.5 }{\em The function $R_g$ always belongs to
$L^1$ and its $L^1$-norm does not exceed $2 ||\eta||$. }

\medskip The classical Boole theorem says that if $d\eta$ 
is non-negative and pure singular, then $N_g(s)=\eta (\R)/s$, and 
therefore $R_g$ vanishes identically. The next two corollaries can 
be viewed as quantitative generalizations of this fact:

\medskip\par\noindent{\bf Corollary 3.6 }{\em If $d\eta\ge 0$, then 
$R_g(t)$ is non-negative as well, and} 
$||R_g||_{L^1} = \eta_{\rm{a.c.}}(\R)$.

\medskip\par\noindent{\bf Corollary 3.7 }{\em If $d\eta$ is pure 
singular, then $R_g(t)$ is non-positive and}  
$||R_g||_{L^1} = ||\eta|| - |\eta (\R)|$.

\medskip
For other recent results obtained with 
the help of the logarithmic determinant we refer to 
\cite{Klemes}, \cite{MS1} and \cite{MS3}.
 
\bigskip\par\noindent{\bf \S 4}

\medskip\par\noindent
In \S 2 we used the subharmonic function technique for proving a theorem 
about the Hilbert transform. The idea of logarithmic determinants also 
provides us with a connection which works in the opposite direction: 
starting with a known result about the Hilbert transform, one arrives at 
a plausible conjecture about a non-negative subharmonic function in $\C$ 
represented by a canonical integral of genus one. For illustration, 
we consider a well known inequality
$$
m_f(\lambda) \lesssim
\frac{1}{\lambda^2}\, \int_0^\lambda sm_g(s) ds
+ \frac{1}{\lambda}\, \int_\lambda^\infty m_g(s) ds\,, 
\qquad 0<\lambda<\infty\,,
\leqno (4.1)
$$
where $f=g+i\h g$, $g$ is a test function on $\R$,  
$m_f(\lambda) = \hbox{meas}\{|f|\ge \lambda\}$, and 
$m_g(\lambda) =  \hbox{meas}\{|g|\ge \lambda\}
= N_g(\la) -N_g(-\la)$. 
Inequality (4.1) contains as special cases Kolmogorov's weak $L^1$-type 
inequality $\la m_f(\la) \lesssim ||g||_{L^1}$, and M.~Riesz' inequality
$||f||_{L^p} \lesssim_p ||g||_{L^p}$, $1<p\le 2$. 
Inequality (4.1) can be justly attributed to Marcinkiewicz.
He formulated his general interpolation theorem for sub-linear operators 
in \cite{M}, the proof was supplied by Zygmund in \cite{Z} with 
reference to Marcinkiewicz' letter. Its main ingredient is a decomposition
$g=g\chi_{\{ |g|<\la\}} + g\chi_{\{|g|\ge \la\}}$, where $\chi_E$ is a 
characteristic function of a set $E$. This decomposition immediately 
proves (4.1), see \cite[Section~V.C.2]{Koosis}.

Define a logarithmic determinant $u_f$ of genus one as in (2.7),
and denote by $d\mu_f$ its Riesz measure (i.e. $1/(2\pi)$ times the 
distributional Laplacian $\Delta u_f$). For each Borelian subset 
$E\subset \C$, $\mu_f(E)=\hbox{meas}(f^{-1}E^*)$, where $E^*=\{z:\, 
z^{-1}\in E\}$, and $f^{-1}E^*$ is the full preimage of $E$ under $f$.  
Now, we can express the RHS and the LHS of inequality (4.1) in the 
terms of $\mu_f$.
First, observe that the counting function of $\mu_f$ equals
$$
\mu_f(r) \stackrel{def}= \mu_f\{|z|\le r \} =
\hbox{meas} \{|f(t)| \ge r^{-1}\} = m_f(r^{-1})\,. 
$$
In order to write down $m_g$ in terms of $\mu_f$, we introduce the 
Levin-Tsuji counting function (cf. \cite{Tsuji}, \cite{GO}):
\begin{eqnarray*}
\n_f(r) &=&
\mu_f \{ |z - ir/2| \le r/2\} +
\mu_f \{ |z + ir/2| \le r/2\} \\
&=& \mu_f\{ |\hbox{Im}(z^{-1})|\ge r^{-1}\}
= \hbox{meas} \{|g|\ge r^{-1}\} = m_g(r^{-1})\,.
\end{eqnarray*}
Now, we can rewrite (4.1) in the form:
$$
\mu_f(r) \lesssim r \int_0^r \frac{\n_f(t)}{t^2}\, dt + r^2 \int_r^\infty
\frac{\n_f(t)}{t^3}\, dt \,,
\qquad 0<r<\infty\,.
\leqno (4.2)
$$
We shall show that (4.2) persists for any 
non-negative in $\C$ subharmonic function represented by a canonical 
integral of genus one. In this case the operator $g\mapsto \h 
g$ disappears, and the Marcinkiewicz argument seems to be unapplicable 
anymore.  

Let
$$
u(z)=\int_{\bf C} K(z/\zeta)\, d\mu(\zeta)\,,
\leqno (4.3)
$$
where $d\mu$ is a non-negative locally finite measure on $\bf C$ such that
$$
\int_{\bf C}
\min\left(\frac{1}{|\zeta|},\frac{1}{|\zeta|^2}\right)\,d\mu(\zeta)<\infty\,.
\leqno (4.4)
$$
Subharmonic functions represented in this form are called {\em canonical 
integrals of genus one}. 
 
Let $M(r,u)=\max_{|z|\le r} u(z)$. 
A standard  estimate of the kernel
$$
K(z) \ll \frac{|z|^2}{1+|z|}\,, \qquad z\in{\bf C}\,,
$$
yields Borel's estimate (cf. \cite[Chapter~II]{GO})
$$
M(r,u) \ll
r\int_0^r \frac{\mu(t)}{t^2}\,dt
+ r^2 \int_r^\infty \frac{\mu(t)}{t^3}\,dt \,.
$$
In particular,
$$
M(r,u) = \left\{
\begin{array}{ll}
o(r), & r\to 0\, \\ \\
o(r^2), & r\to\infty\,.
\end{array}
\right.
$$

\medskip\par\noindent{\bf Theorem~4.5 }{\em 
Let $u(z)\ge 0$ be a canonical integral (4.3) of genus one, then}
$$
M(r,u) \ll r\int_0^r \frac{\n(t)}{t^2}\,dt + r^2 \int_r^\infty
\frac{\n(t)}{t^3}\, dt\,.
\leqno (4.6)
$$

\medskip
The RHS of (4.6) does not depend on the bound for 
the integral (4.4), this makes the result not so obvious. 
By Jensen's formula, $\mu (r) \le M(er,u)$, so that $\mu (r) 
\lesssim \hbox{the RHS of (4.6)}$.
As a corollary we immediately obtain (4.2) and 
the Marcinkiewicz estimate (4.1).

\bigskip\par\noindent{\bf \S 5}
 
\medskip\par\noindent
Here we sketch the proof of Theorem~4.5.

We shall need two auxiliary lemmas. The first one is a version of the 
Levin integral formula without remainder term
(cf. \cite[Section~{\rm IV}.2]{Levin1}, 
\cite[Chapter~1]{GO}). The proof can be found in
\cite{MOS}

\medskip\par\noindent{\bf Lemma~5.1 }{\em Let $v$ be a subharmonic 
function in $\bf C$ such that 
$$
\int_0^{2\pi} |v(re^{i\theta})|\,|\sin\theta|\, d\theta
= o(r)\,, \qquad r\to 0\,,
\leqno (5.2)
$$
and 
$$
\int_0 \frac{\n(t) + v^-(t) + v^-(-t)}{t^2}\, dt < \infty\,.
\leqno (5.3)
$$
Then
$$
\frac{1}{2\pi}
\int_0^{2\pi} v(Re^{i\theta}|\sin\theta|) \frac{d\theta}{R\sin^2\theta}
= \int_0^R \frac{\n(t)}{t^2}\, dt\,, \qquad 0<R<\infty\,,
\leqno (5.4)
$$
where $\n(t)$ is the Levin-Tsuji counting function, and the integral
on the LHS is absolutely convergent.}

\medskip
The next lemma was proved in a slightly different setting in 
\cite[\S 2]{LO}, see also \cite[Lemma~5.2, Chapter~6]{GO}

\medskip\par\noindent{\bf Lemma~5.5 }{\em 
Let $v(z)$ be a subharmonic function in $\bf C$ satisfying conditions 
(5.2) and (5.3) of the previous lemma, let
$$
T(r,v) = \frac{1}{2\pi} \int_0^{2\pi} v^+(re^{i\theta})\,d\theta
$$    
be its Nevanlinna characteristic function, and let
$$
\T(r,v) = \frac{1}{2\pi} \int_0^{2\pi} v^+(re^{i\theta}|\sin\theta)|)
\frac{d\theta}{r\sin^2\theta}
$$
be its Tsuji characteristic function.
Then}
$$
\int_R^\infty \frac{T(r,v)}{r^3}\, dr \le \int_R^\infty
\frac{\T(r,v)}{r^2}\, dr\,, \qquad 0<R<\infty\,.
\leqno (5.6)
$$

\medskip For the reader's convenience, we recall the proof.
Consider the integral
$$
I(R)=\frac{1}{2\pi} \int\!\!\!\int_{\Omega_R}
\frac{v^+(re^{i\theta})}{r^3}\,
dr\,d\theta\,,
$$
where $\Omega_R=\{z=re^{i\theta}:\, r>R|\sin\theta|\}
 = \{z: |z\pm iR/2|>R/2\}$. Introducing a new variable
$\rho=r/|\sin\theta|$ instead of $r$, we get
$$
I(R) = \int_R^\infty \frac{d\rho}{\rho^2} 
\left\{
\frac{1}{2\pi} \int_0^{2\pi} v^+(\rho|\sin\theta|
e^{i\theta}) \frac{d\theta}{\rho\sin^2\theta}
\right\} = \int_R^\infty \frac{\T(\rho, v)}{\rho^2}\, d\rho\,.
$$
Now, consider another integral
$$
J(R) = \frac{1}{2\pi} \int\!\!\!\int_{K_R} \frac{v^+(re^{i\theta})}{r^3}\,
dr\, d\theta\,, 
$$ 
where $K_R=\{z: |z|> R\}$. Since $K_R\subset \Omega_R$, we have $J(R)\le
I(R)$. Taking into account that
$$
J(R)=\int_R^\infty \frac{dr}{r^3} 
\left\{
\frac{1}{2\pi} \int_0^{2\pi} v^+(re^{i\theta}) \,d\theta
\right\} = \int_R^\infty \frac{T(r,v)}{r^3}\, dr
$$
we obtain (5.6). $\Box$

\medskip\par\noindent{\em Proof of Theorem~4.5: }
Due to Borel's estimate condition (5.2) is fulfilled. Due to 
non-negativity of $u$ and (4.4), condition (5.3) holds as well.
Using monotonicity of $T(r,u)$, Lemma~5.5, and then Lemma~5.1, we obtain
\begin{eqnarray*}
\frac{T(R,u)}{R^2} &\le& 2 \int_R^\infty \frac{T(r,u)}{r^3}\, dr \\ \\
&\stackrel{(5.6)}\le& 2
\int_R^\infty \frac{\T(r,u)}{r^2}\, dr \\ \\
&\stackrel{(5.4)}=& 2 \int_R^\infty \frac{dr}{r^2}
\int_0^r
\frac{\n(t)}{t^2}\, dt  \\ \\
&=& \frac{2}{R} \int_0^R \frac{\n(t)}{t^2}\, dt + 2\int_R^\infty
\frac{\n(t)}{t^3}\, dt\,.
\end{eqnarray*}  
The inequality $M(r,u)\le 3T(2r,u)$ completes the
job. $\Box$

\bigskip\par\noindent{\bf \S 6}

\medskip\par\noindent Non-negativity of $u(z)$ in $\C$ seems to be too 
strong assumption, a more natural one is non-negativity of $u(x)$ on $\R$. 

\medskip\par\noindent{\bf Theorem~6.1 }{\em 
Let $u(z)$ be a canonical integral (4.3) of genus one, and
let $u(x)\ge 0$, $x\in \R$. Then
$$
M(r,u) \ll r^2
\left[ \int_r^\infty \frac{\sqrt{\n^*(t)}}{t^2}\,dt \right]^2\,,
\leqno (6.2)
$$
where}
$$
\n^* (r) = r \int_0^r \frac{\n(t)}{t^2}\, dt +
r^2 \int_r^\infty \frac{\n (t)}{t^3}
\left(1+\log\frac{t}{r} \right)\,dt
\leqno (6.3)
$$

\medskip The proof of Theorem~6.1 is given in \cite{MOS}. The method of 
proof differs from that of Theorem~4.5, and is more technical than one 
would wish.

Fix an arbitrary  $\epsilon>0$.
Then by the Cauchy inequality
\begin{eqnarray*}
\left[ \int_r^\infty \frac{\sqrt{\n^*(t)}}{t^2}dt \right]^2
&=&
\left[ \int_r^\infty
\frac{\sqrt{\left(1 +
\log^{1+\epsilon}\frac{t}{r}\right)\n^*(t)}}{t^{3/2}}\,
\frac{dt}{t^{1/2}\sqrt{1+\log^{1+\epsilon}\frac{t}{r}}}
\right]^2 \\ \\
&\ll_{\epsilon}& \int_r^\infty \frac{\n^*(t)}{t^3}
\left(1+\log^{1+\epsilon}\frac{t}{r} \right)\,dt \\ \\
&\ll_{\epsilon}&
\frac{1}{r} \int_0^r \frac{\n(s)}{s^2}\,ds +
\int_r^\infty \frac{\n(s)}{s^3} \left(1+\log^{3+\epsilon}\frac{s}{r}
\right)\, ds\,.
\end{eqnarray*}
Thus we get

\medskip\par\noindent{\bf Corollary 6.4 }{\em For each $\epsilon>0$,
$$
M(r,u) \ll_{\epsilon}
r \int_0^r \frac{\n(t)}{t^2}\,dt + r^2
\int_r^\infty \frac{\n(t)}{t^3} \left(1+\log^{3+\epsilon}\frac{t}{r}
\right)\, dt\,.
\leqno (6.5)
$$
}

\medskip
Estimate (6.5) is slightly weaker than (4.6); however, it suffices for
deriving inequalities of M.~Riesz and Kolmogorov.
Using Jensen's estimate $\mu (r) \le M(er,u)$, we arrive at

\medskip\par\noindent{\bf Corollary~6.6 }{\em The following inequalities 
hold for canonical integrals of genus one which are non-negative on the 
real axis:

\par\noindent M. Riesz-type estimate:
$$
\int_0^\infty \frac{\mu(r)}{r^{p+1}}\, dr
\ll_p\int_0^\infty \frac{\n (r)}{r^{p+1}}\, dr\,,
\qquad 1<p<2\,,
\leqno (6.7)
$$
weak $(p,\infty)$-type estimate:
$$
\sup_{r\in (0,\infty)} \frac{\mu(r)}{r^p} \ll_p
\sup_{r\in (0,\infty)} \frac{\n (r)}{r^p}\,,
\qquad 1<p<2\,,
\leqno (6.8)
$$
and
Kolmogorov-type estimate:}
$$
\sup_{r\in (0,\infty)} \frac{\mu(r)}{r} \ll
\int_0^\infty \frac{\n (r)}{r^2}\, dr\,.
\leqno (6.9)
$$

\medskip\par\noindent{\bf Remark~6.10 } If the integral on the 
RHS of (6.9) is finite, then $u(z)$ has positive harmonic majorants in 
the upper and lower half-planes which can be efficiently estimated near 
the origin and infinity, see \cite[Theorem~3]{MOS}.

\bigskip\par\noindent{\bf \S 7}
\medskip\par\noindent 
Here we mention several questions related to our results.

\medskip\par\noindent{\bf 7.1 } 
How to distinguish the logarithmic determinants (2.7) of 
$f=g+i\h g$ from other canonical integrals (4.3) 
which are non-negative in $\C$? In other words,
let $dm_f$ be a distribution measure of $f$; i.e. a
locally-finite non-negative measure in $\C$ defined by
$m_f(E) = \hbox{meas}\{t\in\R:\, f(t)\in E \}$ for an arbitrary borelian
subset $E\subset \C$. It should to be interesting to find 
properties of $dm_f$ which do not follow only from non-negativity of the
subharmonic function $u_f(z)$. A similar question can be addressed to 
analytic functions $f(z)$ of Smirnov's class in the unit disk.

\medskip\par\noindent{\bf 7.2 } Let $X$ be a rearrangement invariant 
Banach space of measurable functions on $\R$.
That is, the norm in $X$ is the same for
all rearrangements of $|g|$, and $||g_1||_X\le ||g_2||_X$ provided that
$|g_1|\le |g_2|$ everywhere.
For which spaces does the inequality
$$
||\h g_d||_X \le C_X ||\h g||_X
$$
hold? This question is interesting only for spaces $X$ where the Hilbert
transform  is unbounded; i.e. for spaces which are close in a certain
sense either to $L^1$ or to $L^\infty$. Some natural restrictions on 
$X$ can be assumed: the linear span of the characteristic functions 
$\chi_E$ 
of bounded measurable subsets $E$ is dense in $X$, and $||\chi_E||_X\to 
0$, when $\rm{meas}(E)\to 0$, see \cite[Chapter~3]{BS}.

\medskip\par\noindent{\bf 7.3 } We do not know how to extend estimate 
(1.2) (as well as (1.8)) to more general operators like the 
maximal Hilbert transform, the non-tangential maximal conjugate harmonic 
function, or Calderon-Zygmund operators. A similar question can be 
naturally posed for the Riesz transform \cite{Stein}.

\medskip\par\noindent{\bf 7.4 } Does Marcinkiewicz-type inequality (4.6) 
hold under the assumption that a canonical integral $u$ of genus one is 
non-negative on $\R$? 
According to a personal communication from A.~Ph.~Grishin, the exponent
$3+\epsilon$ can be improved in (6.5). However, his technique also does 
not allow to get rid at all of the logarithmic factor.

\medskip\par\noindent{\bf 7.5 }
Let $u(z)$ be a non-negative subharmonic function in $\mathbb C$,
$u(0)=0$. As before, by $\mu (r)$ and $\n(r)$ we denote the conventional
and the Levin-Tsuji counting functions of the Riesz measure $d\mu$ of $u$.
Assume that $\mu (r) = o(r)$, $r\to 0$. This condition is needed to
exclude from consideration the function  $u(z)=|\hbox{Im}(z)|$ which is
non-negative in ${\bf C}$ and harmonic outside of $\bf R$.
Let $\mathcal M$, ${\mathcal M}(0)=0$, ${\mathcal M}(\infty) =
\infty$, be a (regularly growing) majorant for $\n (r)$.
What can be said  about the majorant
for $\mu (r)$? If ${\mathcal M}(r) = r^p$, $1<p<\infty$,
then we know the answer:
$$
\sup_{r\in (0,\infty)} 
\frac{\mu (r)}{r^p}
\le C_p \sup_{r\in (0,\infty)}
\frac{\n(r)}{r^p}\,,
$$
and
$$
\int_0^\infty \frac{\mu (r)}{r^{p+1}}\, dr
\le C_p \int_0^\infty \frac{\n(r)}{r^{p+1}}\, dr\,.
$$
It is more difficult and interesting
to study majorants ${\mathcal M}(r)$ which grow faster than any power of
$r$ when $r\to \infty$,
and decay to zero faster than any power of $r$ when $r\to 0$. The
question might be related   
to the classical Carleman-Levinson-Sjoberg
``$\log\log$-theorem'', and the progress may lead to new results about the 
Hilbert transform.

\bigskip\par\noindent Vladimir Matsaev: {\em School of Mathematical
Sciences,
Tel-Aviv University, \newline
\noindent Ramat-Aviv, 69978, Israel

\par\noindent matsaev@math.tau.ac.il}

\bigskip\par\noindent Iossif Ostrovskii: {\em
Department of Mathematics, Bilkent University, \newline
06533 Bilkent, Ankara, Turkey

\par\noindent iossif@fen.bilkent.edu.tr}

\medskip\par\noindent and {\em Verkin Institute for Low Temperature
Physics and Engineering, \newline
310164 Kharkov, Ukraine

\par\noindent ostrovskii@ilt.kharkov.ua}

\bigskip\par\noindent Mikhail Sodin: {\em School of Mathematical
Sciences, Tel-Aviv University, \newline 
\noindent Ramat-Aviv, 69978, Israel

\par\noindent sodin@math.tau.ac.il}

\end{document}